\newcommand{\Y}{{\Bbb S}^1 \times \Sigma}

\newcommand{\Z}{\Bbb Z}
\newcommand{\R}{\Bbb R}

\newcommand{\D}{\cal D}
\newcommand{\DD}{\Bbb D}
\renewcommand{\S}{\Sigma}
\newcommand{\CP}{\Bbb C \Bbb P}
\newcommand{\SS}{{\Bbb S}^1}

\newcommand{\M}{M_{\Sigma}}

\renewcommand{\a}{\alpha}
\renewcommand{\b}{\beta}
\newcommand{\g}{\gamma}
\newcommand{\x}{\times}

\newcommand{\bd}{\partial}

\newcommand{\Diff}{\operatorname{Diff}}
\newcommand{\isom}{\cong}
\newcommand{\coef}{\operatorname{coef.\ of }}
\newcommand{\ar}{\rightarrow}
\newcommand{\G}{\Gamma}
\newcommand{\pari}[1]{{\partial \over {\partial #1}}}

\documentstyle[12pt,amscd,amssymb]{amsart}

\theoremstyle{plain}
\begingroup
\theorembodyfont{\sl}
\newtheorem{thm}{Theorem}
\newtheorem{cor}[thm]{Corollary}
\newtheorem{lem}[thm]{Lemma}
\newtheorem{prop}[thm]{Proposition}
\newtheorem{conj}[thm]{Conjecture}
\endgroup

\theoremstyle{definition}
\newtheorem{defn}[thm]{Definition}

\theoremstyle{remark}
\newtheorem{rem}[thm]{Remark}

\title[Donaldson invariants for some glued manifolds]
{Donaldson invariants for some glued manifolds}

\author{Vicente Mu\~noz}
\address{Mathematical Institute \\
24-29 St Giles'  \\
Oxford, OX1 3LB}
\email{vmunoz@@maths.ox.ac.uk}
\thanks{\hbox{$^*$}Supported by a grant from Banco de Espa\~na}
\date{October 1995}

\begin{document}

\baselineskip.55cm

\maketitle
\begin{abstract}
  We prove that every suitable $4$-manifold with $b_1=0$ and with an
  embedded Riemann surface of genus $2$ is of simple type. We find a
  relationship between the basic classes of two of these $4$-manifolds
  and those of the connected sum along the Riemann surface.
\end{abstract}

\section{Notations and facts}
\label{sec:notations}

This is a second paper on the problem of studying the behaviour of the
basic classes of a $4$-manifold under a determined surgical
operation. In~\cite{Munoz} we studied how the Seiberg-Witten
invariants behaved under connected sum of two $4$-manifolds along a
Riemann surface. Here we try to find
the Donaldson
invariants (to a reasonable extent) under the same surgery. We also
check the agreement with the conjecture of Witten~\cite{Witten} about
the relation between the Seiberg-Witten and the Donaldson invariants.
The differences in our results appear regarding two aspects:

\begin{itemize}
\item The computations in the instanton case are much more cumbersome
  which hampers us to extend our results beyond the case $g=2$.
\item On the other hand, the gluing theory is more developed in the
  instanton case than in the Seiberg-Witten case, so we can obtain
  information also for basic classes $\kappa$ with $| \kappa \cdot \S
  | < 2g-2$.
\end{itemize}

Before embarking in the statement of the results, we gather briefly
some notation which will be useful. This follows mainly~\cite{KM}.
We say that a four-manifold $X$ is {\em suitable\/} when
$b^+ -b_1$ is odd and $b^+ > 1$. We have defined for any
$w \in H^2(X ;\Z)$ and $p_1 \in H^4(X ;\Z) \isom \Z$ with $p_1 \equiv
w^2 \pmod 4$ the Donaldson invariant $D_X^{w,d}$ as a symmetric
polynomial on ${\Bbb A(X)}= Sym^*(H_{even}(X)) \otimes \bigwedge^*(H_{odd}(X))$
of degree $2d = -2p_1-3(1-b_1+b^+)$ (whenever this quantity
is positive). The degree of the elements in $H_i$ is $4-i$.
We fix $w$. There are several ways of wrapping up all the information
about the
different $p_1$ in a single series, defining so (the class of the
point is $x$):
$$ D_X^w(e^{t\a}) = \sum { D_X^{w,d}(\a ^d) \over d!} t^d $$
$$ \D_X^w(\a)=D_X^w(e^{t\a+\lambda x}) = \sum { D_X^{w,d+2a}(\a ^d x^a)
  \over d! \, a!} t^d\lambda^a $$
for $\a \in H^2(X)$ (our main concern will be the invariants for the
second homology as we will soon suppose $b_1=0$).

We say that $X$ is of {\em $w$-simple type\/} when $x^2 -4$ annihilates the
Donaldson invariant $D_X^w$, that is when $D_X^w((x^2 -4)z)=0$ for
all $z \in {\Bbb A(X)}$. It is a fact that when $X$ is of $w$-simple
type for some $w$, it is so for every $w'$, and it is called of simple
type for brevity. In this case
$$\D_X^w(\a) =D_X^w(e^{t\a})\cosh 2\lambda +D_X^w({x \over
  2}e^{t\a})\sinh 2\lambda$$
Kronheimer and Mrowka defined in~\cite{KM} another series containing the
same information
$$ \DD_X^w(\a)=D_X^w(e^{t\a})+D_X^w({x \over 2}e^{t\a})$$

\begin{thm}
\label{thm:1}
  Let $X$ be a suitable manifold with $b_1=0$ and of simple type. Then
  we have
  $$\DD_X^w(\a)= e^{Q(\a) /2} \sum a_{i,w} e^{K_i \cdot \a}$$
  for finitely many cohomology classes $K_i$ (called {\em basic
  classes\/}) and rational numbers $a_{i,w}$. These classes are lifts
  to integral homology of $w_2(X)$. Moreover, for any embedded surface
  $ S \hookrightarrow X$ of genus $g$ and positive self-intersection,
  one has $2g-2 \geq S^2 +|K_i \cdot S|$. Also the rational numbers
  satisfy $a_{i,w}=(-1)^{{K_i \cdot w +w^2} \over 2} a_i$, for some $a_i \in
  {\Bbb Q}$.
\end{thm}

\begin{defn}
  Here we are going to have a situation with $w, \S \in H^2(X;
  \Z)$ satisfying $w \cdot \S \equiv 1 \pmod 2$ and $\S^2 =0$. Then we
  define $$D_X=D^w_X +D^{w+\S}_X$$
  Obviously, $D_X$ depends on $w \mod \S$.
  Note that since $(w+\S)^2 \equiv w^2 +2 \pmod 4$, we can recuperate
  $D^w_X$ and $D^{w+\S}_X$ from the series $D_X$. Analogously we
  define $\D_X =\D^w_X + \D^{w+\S}_X$.
\end{defn}

\begin{prop}
  Suppose $X$ is a suitable manifold of simple type with $b_1=0$ and $\DD_X^w =
  e^{Q/2} \sum a_i e^{K_i}$. Then setting $d_0= -w^2 -{3 \over
    2}(1+b^+)$ we have $$D_X(e^{\a})=e^{Q({\a})/2} \sum_{K_i \cdot \S
    \equiv 2 \pmod 4} a_i e^{K_i \cdot \a} + e^{Q/2} \sum_{K_i \cdot
    \S \equiv 0 \pmod 4} i^{-d_0} a_i e^{iK_i \cdot \a}$$
  So it is equivalent giving $\DD_X^w $ to giving $D_X$.
\end{prop}

\begin{pf}
  Since $((w+\S)^2+K_i \cdot (w+\S)) -(w^2 + K_i \cdot w)) = 2(w \cdot
  \S + K_i \cdot \S /2)$ we have
  $$\DD_X^{w+\S} = e^{Q/2} \sum_{K_i \cdot \S
    \equiv 2 \pmod 4} a_i e^{K_i} -
  e^{Q/2} \sum_{K_i \cdot \S
    \equiv 0 \pmod 4} a_i e^{K_i}$$
  Now $D_X^w(e^{\a})= {1 \over 2}(\DD_X^w(\a) + i ^{-d_0}\DD_X^w(i\a))$
  and $D_X^{w+\S}(e^{\a})= {1 \over 2}(\DD_X^{w+\S} (\a) - i
  ^{-d_0}\DD_X^{w+\S} (i\a))$, so
  $$D_X(e^{\a})=e^{Q({\a})/2} \sum_{K_i \cdot \S
    \equiv 2 \pmod 4} a_i e^{K_i \cdot \a} + i^{-d_0}e^{Q/2} \sum_{K_i \cdot
    \S \equiv 0 \pmod 4} a_i e^{iK_i \cdot \a}$$
\end{pf}

\section{Statement of results}
\label{sec:results}

There are many ways of cooking different four-manifolds out of old
ones. It is an interesting problem to relate the invariants of the new
manifolds with those of the given ones (and to find the basic classes
when appropriate) (for very nice examples see~\cite{Stipsicz}).

The case this paper is concerned with is the following. Let $\bar X_1$
and $\bar X_2$ be two suitable manifolds with $\S_i \subset \bar X_i$
embedded surfaces of the same genus $g$ and self-intersection zero. If we
remove a tubular neighbourhood of $\S_i$ we end up with a
manifold $X_i$ with boundary $Y= \Y$. Choosing an identification
between the boundaries of $X_1$ and $X_2$ (there is an element of
choice here which produces different four-manifolds~\cite{Gompf}), we
can construct another manifold $X=X_1 \cup_Y X_2$.
The (closure of the) tubular neighbourhood removed is always
diffeomorphic to $A= D^2 \x \S$. Obviously $\bar X_i = X_i \cup_Y A$.

Our starting hypotheses are that $b_1=0$ for both $\bar X_i$ and that
there exists $w_i \in H^2(\bar X_i; \Z)$ with $w_i \cdot \S_i \equiv 1 \pmod
2$, equivalently the cohomology class of $\S_i$ is an odd multiple of a
non-torsion primitive class. We fix $w$'s in all the manifolds involved
($\bar X_1$, $\bar X_2$ and $X$)
once and for all in a compatible way (i.e. the restriction of $w$ to
$X_i$ coincides with the restriction of $w_i$ to $X_i$) with
$w^2=w_1^2+w_2^2 \pmod 4$ and
such that $ w_i \cdot \S_i \equiv 1 \pmod 2$. We drop the subindices, so
we will not differentiate
these $w$'s, since the context makes always clear to which manifold they
refer.

The series $\D_X^w(\a)$ is determined by its action on elements $\a$ lying
in:
\begin{itemize}
\item $H^2(X_1) \oplus H^2(X_2) \subset H^2(X)$.
\item elements $D$ such that $D=D_1 +D_2$ with $\bd D_1 = -\bd D_2$ being
  $\SS \subset Y=\Y$ (or a multiple of $\SS$).
\item elements $D_{\gamma}$ such that $D_{\gamma}=D_1 +D_2$ with $\bd
  D_1 = -\bd D_2$ representing a cycle $\gamma$ in the subgroup $
  H_1( \S; \Z) \subset H_1( \Y; \Z)$.
\end{itemize}

For studying the behaviour on elements of the first kind we can
suppose that $\a \in H^2(X_1)$ without loss of generality. As
  explained in~\cite{Braam}~\cite{Don2} we
  have Floer homology groups $HF_*(Y)$, cooked out of the flat connections on
  $Y$. The relevant issues we must be careful with are the appearance of
reducibles on $X_i$ (but in our case they are ruled out since the moduli space
$\M$ of flat connections on $\S$ consists entirely of irreducibles)
and of flat connections on $X_i$ (they present no problem using the
blow-up trick of~\cite{MM}).
One has an invariant $D^w_{X_1}(\a^d)
  \in HF_*(Y)$ and another invariant $D^w_{X_2}(1) \in HF^*(Y)$ such
  that $D_X(\a^d)=<D^w_{X_1}(\a^d), D^w_{X_2}(1)>$. Looking at it more
  closely we realize that the pairing above gives either $D^w_X(\a^d)$ or
  $D^{w+\S}_X(\a^d)$ depending on $d \mod 4$.

For dealing with the other two cases there is an analogous approach
(see~\cite{Braam}). Let $D=D_1 +D_2$ with $\bd
  D_1 = -\bd D_2 = \gamma \in  H_1( \Y; \Z)$. Then we have
  Fukaya-Floer groups $HFF_*(Y, \gamma)$ such that the invariants are
$D^w_{X_1}(D_1^i) \in HFF_*(Y, \gamma)$ and
$D^w_{X_2}(D_2^j) \in HFF^*(Y, \gamma)$ with
$$D_X (D^m)=\sum
{\left(\begin{array}{c} m \\ i \end{array} \right)}
 <D^w_{X_1}(D^i), D^w_{X_2}(D^{m-i})>$$
We aim to prove that:

\begin{thm}
  Let $X$ have $b_1=0$ and an embedded surface of genus $2$ and
  self-intersection zero. Then $X$ is of
  simple type.
\end{thm}

So let $\bar X_i$ be two suitable manifolds as above with
embedded surfaces $\S_i \subset \bar X_i$ of genus
two and self-intersection zero (and being an odd multiple of a
non-torsion class in homology). Form $X= \bar X_1 \#_{\S} \bar X_2$
removing tubular neighbourhoods of $\S_i$ and identifying the
boundaries $Y=\Y$ (recall that the
diffeomorphism type of the resulting manifold depends on the homotopy
class of the gluing identification). The manifold $X$ is of simple type
by the last theorem.

\begin{thm}
\label{thm:main}
  Suppose that $\bar X_i$ have $b_1=0$ and are of simple type.
  Then $X$ is of
  simple type and every basic class $\kappa$ of $X$ intersects $Y$ in
  $n\SS$ where $n=0,2,-2$. Moreover the sum of the coefficients
  $c_{\kappa}$ of the different basic classes $\kappa$ agreeing with
  $(K,L) \in H^2(\bar X_1;\Z)/\Z[\S_1] \oplus H^2(\bar
  X_2;\Z)/\Z[\S_2]$ is zero unless $(K,L)$ comes from $(K_i,L_j) \in
  H^2(\bar X_1;\Z) \oplus H^2(\bar X_2;\Z)$ where $K_i$ is a basic
  class for $\bar
  X_1$, $L_j$ is a basic class for $\bar X_2$ and $K_1 \cdot \S_1 =
  L_j \cdot \S_2 = \pm 2$, in which case it is $\pm 32$ times the
  product of the coefficients of $K_i$ and $L_j$.
\end{thm}

For $g>2$ we have the obvious following

\begin{conj}
  \label{conj}
   Let $\bar X_i$ have $b_1=0$ and be of simple type.
   Suppose that there are embedded surfaces $\S_i \subset \bar X_i$ of genus
  $g$ and self-intersection zero. Form $X= \bar X_1 \#_{\S} \bar X_2$.
  Then $X$ is of
  simple type and every basic class $\kappa$ of $X$ intersects $Y$ in
  $n\SS$ where $n$ is an even integer with $-(2g-2) \leq n \leq
  (2g-2)$.
  Moreover the sum of the coefficients
  $c_{\kappa}$ of the different basic classes $\kappa$ agreeing with
  $(K,L) \in H^2(\bar X_1;\Z)/\Z[\S_1] \oplus H^2(\bar
  X_2;\Z)/\Z[\S_2]$ is zero unless $(K,L)$ comes from $(K_i,L_j) \in
  H^2(\bar X_1;\Z) \oplus H^2(\bar X_2;\Z)$ where $K_i$ is a basic
  class for $\bar
  X_1$, $L_j$ is a basic class for $\bar X_2$ and $K_1 \cdot \S_1 =
  L_j \cdot \S_2 = \pm (2g-2)$, in which case it is $\pm 2^{7g-9}$ times the
  product of the coefficients of $K_i$ and $L_j$.
\end{conj}

Here we would like to point out the close relationship between this
result and the results in~\cite{Munoz} about Seiberg-Witten
invariants. Witten~\cite{Witten} has conjectured that for a simply
connected manifold the condition of being simple
type and Seiberg-Witten simple type are equivalent, that in that
case the basic classes are the same as the Seiberg-Witten
basic classes and that the shape of the Donaldson series is
$$ \DD_X^w = e^{Q/2} \sum a_i e^{K_i}$$
where $$a_i = (-1)^{{K_i \cdot w +w^2} \over 2}2 ^{2+{1 \over 4}(7
  \chi +11 \sigma)} SW_X(K_i)$$

We proved in~\cite{Munoz} that for basic classes $K_i$ such that
$K_1 \cdot \S = K_2 \cdot \S =
\pm (2g-2)$ one has
$$ \sum_{L \in \pi^{-1}(K_1,K_2)} SW_X(L) = SW_{\bar X_1} (K_1)
SW_{\bar X_1} (K_2) $$
Now the topological numbers are (see~\cite{Munoz}) $ \chi_X=
\chi_{\bar X_1} + \chi_{\bar X_2} +4g-4$ and $ \sigma_X =
\sigma_{\bar X_1}+\sigma_{\bar X_2}$, so for $g=2$, ${2+{1 \over 4}(7
  \chi +11 \sigma)}$ is increased in $5$. Also from equation~\eqref{G} any
basic class $K \in \pi^{-1}(K_1,K_2)$ with $K \cdot \S= \pm (2g-2)$ has
$K   \cdot w = K_1 \cdot w_1 + K_2 \cdot w_2 \pm 2 \pmod 4$
(and we chose $w^2=w_1^2+w_2^2 \pmod 4$).
Therefore the sum of the coefficients
of the different basic classes $K$ agreeing with
  $(K_1,K_2) \in H^2(\bar X_1;\Z)/\Z[\S_1] \oplus H^2(\bar
  X_2;\Z)/\Z[\S_2]$ is $\pm 32$ times the
  product of the coefficients of $K_i$ and $L_j$

Note that this result is more restricted in the sense that it is only
valid for the case $g=2$ but it is more general in the sense that it
also gives information about basic classes $K$ with $|K \cdot \S | <
2g-2$.

{\em Acknowledgements:\/} I would like to thank my supervisor Prof.\
Donaldson for
his encouragement and many useful ideas. I am also very grateful to
Banco de Espa\~na for financial support during the year 94/95.

\section{The case $g=2$}
\label{sec:g=2}

We are specially interested in the case in which the genus is $2$
since the computations can then be carried out quite explicitly. For
this we recall the Floer homology of $\Y$ when $g=2$. In general,
there is the isomorphism $ HF^*(\Y) \isom QH^*(\M) $ with the quantum
cohomology of $\M$, the moduli space of rank-$2$ stable bundles over
$\S$ with odd determinant. The universal bundle yields a map
$\mu: H_*(\S) \ar H^*(\M)$
given by slanting with the first Pontrjagin class.
$\M$ can be described
algebraically as the intersection of two quadrics in $\CP^5$. From
this description, Donaldson~\cite{Don1} found the space of lines in
$\M$, the necessary input for finding the quantum corrections.
$QH^*(\M)$ has (integral) homology equal to $\Z$ in degrees $0$,
$2$, $4$ and $6$. The generators are $1$, $h$, $l$ and $p$ which
correspond in the description above to the fundamental class, a
plane, a line and a point in $\M$.
The map $\mu$ gives an isomorphism
$\mu: H_1(\S) \ar H^3(\M)$, describing the other non-zero bit of
the homology of $\M$. We have
$$ h * h = 4l+4 $$
$$ h*h *h = 4p +12 h $$
$$ h*h *h * h =h * h$$

In general we will drop the $*$ symbol for denoting the quantum
product. $QH^*(\M)$ has a natural $\Z/4$-grading coming from
reducing the $\Z$-grading above.
The standard basis is given by $e_i=h^i$, $0 \leq i \leq 3$, and
elements $\mu(\a)$, where $\a$ runs through a basis of $H_1(\S)$. Note
that the matrix $< e_i, e_j>$ is
$$ 4{\left(\begin{array}{cccc} 0&0&0&1 \\ 0&0&1&0 \\ 0&1&0&16 \\ 1&0&16&0
    \end{array} \right)} $$

The pairings $<\mu(\a), e_i>$ are all zero, so $QH^3(\M)$ is
orthogonal to the ``even'' part of $QH^*(\M)$. This is an important
remark since it will allow us to ignore the ``odd'' part in later
computations.

\section{The first case}
\label{sec:1st}

Here we suppose $g \geq 2$.
First we deal with the invariants for classes $\a \in H_2(X_1)
\subset H_2(X)$ and $\b \in H_2(X_2) \subset H_2(X)$.
We fix elements $z_l= \S^a x^b \a_1 \cdots \a_r \in {\Bbb A}(\S)$,
$\a_i \in H_1(\S)$, such
that the corresponding $e_l = \mu(\S)^a \mu(x)^b \mu(\a_1) \cdots
\mu(\a_r)$ form a basis for $HF_*(\M)$ (quantum multiplication is
understood throughout). Then we have the following:
\begin{equation}
  D^w_A(z_l)= e_l
  \label{A}
\end{equation}

For the open manifold $X_1$ we write $D^w_{X_1}(z)=<D^w_{X_1}(z), e_l> e_l^*$,
where $\{e_l^*\}$ is the dual basis of $\{ e_l \}$. Transferring the
cycles contained in $z_l$ to $X_1$ we get for $z  \in {\Bbb A}(X_1)$.
$$D^w_{X_1}(z)=<D^w_{X_1}(z),  D^w_A(z_l)> e_l^* =<D^w_{X_1}(z_lz),
D^w_A(1)> e_l^* =D_{\bar X_1}(z_lz) e_l^*$$

Then
$$D_X(e^{t(\a+\b)})=<D^w_{X_1}(e^{t\a}),D^w_{X_2}(e^{t\b})>=D_{\bar
  X_1}(e^{t\a}z_l) D_{\bar
  X_2}(e^{t\b}z_m) <e_l^*, e_m^*>$$
Therefore the series for $X$ is
determined by the series for both sides.
We can work out $\D_X(\a+\b)$ adding the class of the point to
either side.
Note that on cycles of
the type $T_{\g}= \SS \x \g \subset \Y$ the series is
constant, which agrees with the case of simple type, in which tori of
self-intersection zero have intersection zero with all basic classes
(and so the invariants are zero on such tori).

When both of $\bar X_i$ are of simple type (and $b_1=0$),
we have that $X$ is of
simple type (on cycles of the first kind) and writing
$$\DD^w_{\bar X_1} (e^{\a}) = e^{Q(\a)/2} \sum a_i e^{K_i \cdot \a},
\; \DD^w_{\bar X_2} (e^{\b}) = e^{Q(\b)/2} \sum b_j e^{L_j \cdot \b}$$
(we drop the dependence of $w$ on the coefficients $a_i$ and $b_j$), we
have that
$\DD^w_{\bar X_1} (e^{\a}\S^a) = e^{Q(\a)/2} \sum a_i (K_i \cdot \S)^a
e^{K_i \cdot \a}$. Also
$$D_{\bar X_1} (e^{\a}) = e^{Q(\a)/2} \sum a'_i e^{K'_i \cdot \a}, \;
D_{\bar X_2} (e^{\b}) = e^{Q(\b)/2} \sum b'_j e^{L'_j \cdot \b}$$
where $a'_i=a_i$ and $K'_i =K_i$ when $K_i \cdot \S \equiv 2 \pmod 4$ and
$a'_i=i^{-d_0}a_i$ and $K_i =iK_i$ when $K_i \cdot \S \equiv 0 \pmod 4$ (and
analogously for $\bar X_2$). Now
$$D_X(e^{\a+\b}) = e^{Q(\a+\b)/2}\sum_{l,m, r=r'=0}
a'_ib'_j e^{K'_i \cdot \a +  L'_j \cdot \b}
((K'_i \cdot \S)^a 2^b (L'_j \cdot \S)^{a'}2^{b'} <e_l^*, e_m^*>)$$
$$\DD^w_X(e^{\a+\b}) = e^{Q(\a+\b)/2}\sum_{l,m, r=r'=0}
(a'_ib'_j)' e^{K_i \cdot \a +  L_j \cdot \b}
((K'_i \cdot \S)^a 2^b (L'_j \cdot \S)^{a'}2^{b'} <e_l^*, e_m^*>)$$
here $a$, $b$ and $r$ correspond to $l$ and
$a'$, $b'$  and $r'$ correspond to $m$.
Using $\a=\S$, $\b=0$ and $\a=0$, $\b=\S$ we see that $K_i \cdot \S =
L_j \cdot \S$. Moreover from theorem~\ref{thm:1} this number is even
and less or equal than $2g-2$ in absolute value.

\subsection*{The genus $2$ case}

When $g=2$, $K_i \cdot \S$ has to be $-2$, $0$ or $2$. For the case
$K_i \cdot \S =0$ the number in the right hand side vanishes so we
have:

\begin{thm}
\label{thm:2}
  $$\DD^w_X(e^{\a+\b}) = e^{Q(\a+\b)/2}(\sum_{K_i \cdot \S=L_j \cdot
    \S = 2} 32a_ib_j \, e^{K_i \cdot \a +  L_j \cdot \b}+\sum_{K_i \cdot
    \S=L_j \cdot \S= -2} -32a_ib_j \, e^{K_i \cdot \a +  L_j \cdot \b})$$
\end{thm}

\begin{pf}
  We use the standard basis for $QH^*(\M)$ and note that $h$
  corresponds to $2 \S$ by lemma~\ref{lem:ST}. So in the above expression
  $b=b'=0$ and $0 \leq a,a' \leq 3$. The matrix $<e_l^*, e_m^*>$ is
  $$ {1 \over 4}{\left(\begin{array}{cccc} 0&-16&0&1 \\ -16&0&1&0 \\
    0&1&0&0 \\ 1&0&0&0
    \end{array} \right)} $$
  so when $K_i \cdot \S=2$, the coefficient is computed to be $32$ and
  when $K_i
  \cdot \S=- 2$ it will be $-32$.
\end{pf}

\begin{rem}
  The reason for the signs is easy to work out. First, $w^2$ for $X$
  is congruent $\pmod 2$ with the sum of both of $w^2$ for $\bar X_i$.
  Also $b^+(X) =b^+(\bar X_1)+b^+(\bar X_2)+ (2g-1)$, so
  $-{3 \over 2}(1+b^+(X))= -{3 \over 2}(1+b^+(\bar X_1))-{3 \over
    2}(1+b^+(\bar X_2)) -3(g-1)$.
  Recalling
  that $d_0 = -w^2 -{3 \over 2}(1+b^+)$ and $g=2$,
  we have $d_0(X) \equiv d_0(\bar
  X_1)+d_0(\bar X_2)+1 \pmod 2$. Now the sign comes from the fact that
  the coefficient for the basic class $-K_i$ is $(-1)^{d_0}a_i$, being
  $a_i$ the coefficient for the basic class $K_i$.
\end{rem}

\section{The second case}
\label{sec:2nd}

Let now have a homology class $D$ such that $D=D_1 + D_2$ with $D_i
\in H_2(X_i, \bd X_i; \Z)$ and $\bd D_1 = - \bd D_2 = \SS \subset \Y$
(or a multiple of $\SS$, e.g. when $\S$ is non-primitive).
We can suppose $D^2=0$ (by adding a suitable multiple of $\S$ to $D$).
This will do no harm to the argument and is notationally convenient.
In this section we need to work with the Fukaya-Floer homology version of
the gluing theory. This appears as the limit of a spectral sequence
whose $E_3$ term is $HF_*(Y) \otimes H_*( \CP^{\infty})$. Since all
maps in this spectral sequence
are $H_{odd}(\M) \ar H_{even}(\M)$ and $H_{even}(\M)
\ar H_{odd}(\M)$ and all commute with the action of $\Diff (\S)$
(because the boundary cycle is $\SS$ and therefore invariant under
that group),
they are zero. Thus the spectral sequence degenerates in the
third term and
$$ D^w_{X_i} (D) =( \phi_0, \phi_1, \phi_2, \ldots) \in HF_*(Y)
\otimes H_*( \CP^{\infty}) $$
where we can interpret $\phi_k = D^w_{X_i} (D^k) \in QH_*(\M)$, which
appears in a similar fashion  as the invariant for the closed manifold
but using
the moduli space of connections for the open manifold
with a cylindrical end.
The pairing formula reads
$$D_X (D^m)=\sum
{\left(\begin{array}{c} m \\ i \end{array} \right)}
 <D^w_{X_1}(D^i), D^w_{X_2}(D^{m-i})>$$

Now we use the trick of
 transferring $\S$ from $X_1$ to $X_2$.
$$D_X (D^m \S)=\sum
{\left(\begin{array}{c} m \\ i \end{array} \right)}
 <D^w_{X_1}(D^i \S), D^w_{X_2}(D^{m-i})> = <D^w_{X_1}(D^i), P_{\S}
 D^w_{X_2}(D^{m-i})> $$
for some symmetric map $P_{\S}: H_*(\M) \ar
H_*(\M) $. We have the following result about the structure of
$P_{\S}$ which we do not prove here but hope to return to it later. In
any case, one can avoid using it (see remark~\ref{rem:tr}).

\begin{lem}
\label{lem:PS}
  $P_{\S}$ is quantum multiplication by $\mu(\S)$ plus $i (D
  \cdot \S)$ times the
  identity, i.e. $D^w_{X_1}(D^i\S) = \mu(\S) * D^w_{X_1}(D^i) + i (D
  \cdot \S) D^w_{X_1}(D^{i-1})$.
\end{lem}

\begin{cor}
\label{cor:QS}
  As a consequence of the above lemma we get
  $$D_X(e^{s\S +tD})= <D^w_{X_1} (e^{tD_1} e^{s\S}), D^w_{X_2}(e^{tD_2})> =
    <D^w_{X_1} (e^{tD_1} ), e^{(s\S)(tD_1)} e^{s \mu(\S)} *
    D^w_{X_2}(e^{tD_2})>$$
\end{cor}

\subsection*{The genus $2$ case}

First we start off with a simple result

\begin{lem}
\label{lem:ST}
  For $g=2$ one has $\mu(\S)={1 \over 2}h$ and $\mu(x)=-4h^2 +2$. In
  particular $\mu(x)^2-4=0$.
\end{lem}

\begin{pf}
  Take $X$ to be a K3 surface blown-up in two points. Consider a tight
  surface $S$ of self-intersection $2$ (and therefore of genus $2$) in the
  K3 surface (which existence
  is guaranteed by~\cite{KM}). Let $E_1$ and $E_2$ be the exceptional
  divisors in $X$ and let $\S$ be the proper transform of the tight
  surface, i.e. $\S = S -E_1 -E_2$.
  Put $w=E_1$, so $w \cdot \S =1$ and $\S$ has genus $2$ and
  self-intersection zero. Then $X$ is of simple type,
  $\DD_X^w=e^{Q/2}\cosh E_2 \sinh E_1$. So $\DD_X^w(t\S)=\cosh t \sinh
  t$. The moduli spaces of connections on $X$ are of dimensions $2d
  \equiv 6 \pmod 8$. From all of this we get
  $D^w_X(\S^{3+4n})=2^{2+4n}$.
  Write $X=X_1 \cup_Y A$. Then $\mu(\S)$ is a multiple of $h$, say
  $ah$. Now $D^w_X(\S^{n+6})= <D^w_{X_1}(\S^n) , (ah)^6> =
  <D^w_{X_1}(\S^n) , a^616^2h^2>=a^416^2D^w_X(\S^{n+2})$ (for $n \equiv
  1 \pmod 4$), from where
  $a= {1 \over 2}$.

  For computing $\mu(x)$ we put $\mu(x) = ah^2+b$. We have
  $0=D^w_X((x^{2n}-2^{2n})\S^3) = <D^w_{X_1}(\S^3), (ah^2+b)^{2n}
  -2^{2n}>$, for $n >0$. The first factor in the pairing is of the form
  $ch+d \neq
  0$, since the expression is non-zero for $n=0$. Therefore all the
  second factors (for different $n$) are proportional. This is
  impossible if $b \neq 2$ since $(ah^2+b)^{2n}
  -2^{2n} = {1 \over 16}((16a+b)^{2n} - b^{2n}) h^2 + (b^{2n}-2^{2n})$.
  So $b=2$, $\mu(x)=ah^2+2$ and $\mu(x)^2 -4 = (16a^2 +4a) h^2$
  implying that $0= D^w_X((x^2-4)z)= (16a^2 +4a)D^w_X(\S^2z)$. So
  either $a=0$ or $a=-4$. The first case is impossible since $\mu(x)$
  has a non-trivial component in $H^4$.
\end{pf}

\begin{cor}
  Let $\bar X_1$ have $b_1=0$ and an embedded surface of genus $2$ and
  self-intersection zero (as supposed so far). Then $\bar X_1$ is of
  simple type. \; $\Box$
\end{cor}

Easily we obtain
\begin{equation}
  e^{s\mu(\S)} = 1 + {s \over 2} h + {\cosh 2s -1 \over 16} h^2 + {\sinh
  2s -2s \over 64} h^3
\label{B}
\end{equation}

Let us write
$D^w_{X_1} (e^{tD_1} ) =(\phi_0, \phi_1, \phi_2, \phi_3)$,
$D^w_{X_2} (e^{tD_2} ) =(\psi_0, \psi_1, \psi_2, \psi_3)$ and
$D^w_{A} (e^{t\Delta} ) =(a_0, a_1, a_2, a_3)$
with respect to the standard basis of section~\ref{sec:g=2}
($\Delta =D^2 \x pt \subset D^2 \x
\S =A$). We do not consider the odd part of the Floer homology since
when one of the manifolds (say $\bar X_1$) has $b_1=0$, one has
$<D^w_{X_1} (e^{tD_1} ), (\mu(\a))^*> = 0$ for $\a \in H_1(\S)$, as the
dual of $\mu(\a)$ is of the form $\mu(\b)$. We suppose
that $D \cdot \S =1$
although the general case $D \cdot \S  \neq 0$ is much the same as
this one. We have
$$D_X(e^{s\S +tD}) = $$
$$ = e^{ts}
  (\phi_0 , \phi_1 , \phi_2 ,\phi_3)
  {1 \over 4} {\left(\begin{array}{cccc} 0 & -16 &0 &1 \\  -16 &0 &1
    &0 \\ 0 &1 &0 & 0 \\ 1 &0 & 0 &0 \end{array} \right)}
  {\left(\begin{array}{cccc} 1 & {1 \over 2}s & {1 \over 16}(\cosh 2s
    -1)  &{1 \over 64} (\sinh
    2s -2s)  \\  0 &1 & {1 \over 4} \sinh 2s & {1 \over
    16}(\cosh2s-1) \\ 0 & 0 & \cosh 2s &{1 \over 4}\sinh2s \\ 0 & 0 &
  4 \sinh 2s &
    \cosh 2s \end{array} \right)}
  {\left(\begin{array}{c} \psi_0 \\ \psi_1 \\ \psi_2 \\ \psi_3
    \end{array} \right)}  \\ $$
$$ = e^{ts}
  (\phi_0 , \phi_1 , \phi_2 ,\phi_3)
  {1 \over 4}
  {\left(\begin{array}{cccc}   0 & -16 &0 &1 \\
    -16 & -8s & 1& {1 \over 2}s \\ 0 &1 & {1 \over 4}\sinh 2s & {1
      \over 16}(\cosh 2s -1) \\ 1 &{1 \over 2}s & {1 \over 16}(\cosh 2s
    -1) & {1 \over 64} (\sinh 2s- 2s)  \end{array} \right)}
  {\left(\begin{array}{c} \psi_0 \\ \psi_1 \\ \psi_2 \\ \psi_3
    \end{array} \right)} $$

Call the square matrix in the middle $B$. Now
we can separate according to coefficients corresponding to functions
on $s$.
\begin{equation}
  {\left(\begin{array}{l} \coef e^{2s}e^{ts} \\ \coef e^{-2s}e^{ts} \\
  \coef e^{ts} \\ \coef se^{ts}    \end{array} \right)}
  = {1 \over 4}{\left(\begin{array}{cccc}   0 & 0 & {1 \over 128}
    (4\psi_3+16\psi_2) & {1 \over 128} (4\psi_2+\psi_3) \\0 & 0 &{1 \over 128}
    (4\psi_3- 16\psi_2) & {1 \over 128} (4\psi_2 -\psi_3) \\
    \psi_3-16\psi_1 & \psi_2-16\psi_0 & \psi_1 - {1 \over 16}\psi_3 &  \psi_0 -
    {1 \over 16}\psi_2 \\ 0 & \psi_3-16\psi_1 & 0 & \psi_1 - {1 \over 16}\psi_3
  \end{array} \right)}
  {\left(\begin{array}{c} \phi_0 \\ \phi_1 \\ \phi_2 \\ \phi_3
    \end{array} \right)}
  \label{C}
\end{equation}

Calling the new matrix in the middle $A_{\psi}$, one has
$(e^{2s}, e^{-2s}, 1 ,s)A _{\psi}= (\psi_0, \psi_1, \psi_2 ,\psi_3) B$.
We call $A$ to the matrix $A_a$ corresponding to $A=D^2 \x \S$.

$$  {1 \over 4}{\left(\begin{array}{cccc}   0 & 0 & {1 \over 128}
    (4a_3+16a_2) & {1 \over 128} (4a_2+a_3) \\0 & 0 &{1 \over 128}
    (4a_3- 16a_2) & {1 \over 128} (4a_2 -a_3) \\
    a_3-16a_1 & a_2-16a_0 & a_1 - {1 \over 16}a_3 &  a_0 -
    {1 \over 16}a_2 \\ 0 & a_3-16a_1 & 0 & a_1 - {1 \over 16} a_3
  \end{array} \right)}$$

\begin{lem}
  The matrix $A$ is invertible.
\end{lem}

\begin{pf}
  That the determinant vanishes would imply that either $a_3 = 4a_2$,
  $a_3 = -4a_2$ or $a_3=16a_1$. The first two cases give
  that the first or second row of $A$ is zero respectively, which is
  contradictory as there are examples where the left hand side
  of~\eqref{C} has non-zero two first entries (see remark~\ref{rem:tr}).
  The case $a_3=16a_1$ implies that the series for any
  such $X$ is always of the form
  $e^{ts}(f_1(t)e^{2s}+f_2(t)e^{-2s}+f_3(t))$. This is also valid for
  $X= \CP^1 \x \S =A \cup_Y A$ (see remark~\ref{rem:tr}).
  Particularizing for $t=0$,
  $D_X(e^{s\S})$ is a linear combination of $e^{2s}$, $e^{-2s}$
  and $1$. But from~\eqref{A} and~\eqref{B} we get that
  $$ D_X(e^{s\S}) ={1 \over 16} (\sinh 2s -2s )$$
\end{pf}

As a corollary, the $\phi_i$ are determined by the series $D_{\bar
  X_1}(e^{t{\bar D_1}})$ ($\bar D_1 = D_1 \cup \Delta$). We have
$$D_X(e^{tD})   = <(\phi_0, \phi_1, \phi_2, \phi_3), (\psi_0, \psi_1,
  \psi_2, \psi_3)>= $$
\begin{equation}
  (v_1)^T {\left(\begin{array}{cccc} {32 \over (a_3+4a_2)^2} &0&0&0 \\ 0&
    -{32 \over (a_3-4a_2)^2} &0&0 \\ 0&0&0&{-4 \over (a_3-16a_1)^2} \\
    0&0&{-4 \over (a_3-16a_1)^2} & 8{a_2 -16a_0 \over (a_3-16a_1)^3}
  \end{array} \right)} \, (v_2)
    \label{D}
\end{equation}
where
$$v_i ={\left(\begin{array}{l} \coef e^{2s}e^{ts} \\ \coef e^{-2s}e^{ts} \\
  \coef e^{ts} \\ \coef se^{ts}    \end{array} \right)} $$
for the manifold $\bar X_i$. When $\bar X_1$ is of simple type we can
use $\DD^w_{\bar X_1}$ instead of $D_{\bar X_1}$.
We write $\bar D_i = D_i \cup \Delta$. We have some freedom,
so we impose $ \bar D_i^2 =0$. Then we can write
$$v_1={\left(\begin{array}{c} \sum\limits_{K_i \cdot \S=2} a_ie^{tK_i
      \cdot \bar D_1} \\
  \sum\limits_{K_i \cdot \S=-2}  a_ie^{tK_i \cdot  \bar D_1}  \\
  \sum\limits_{K_i \cdot \S=0}
  a_ie^{tK_i \cdot  \bar D_1}  \\ 0
    \end{array} \right)}$$
When both of $\bar X_i$ are of simple type, we have
\begin{equation}
  \DD^w_X(tD)=
  \label{E}
\end{equation}
$$=(\sum\limits_{K_i \cdot \S=2} \! \! a_ie^{tK_i
      \cdot \bar D_1}, \!
  \sum\limits_{K_i \cdot \S=-2} \!\! \! a_ie^{tK_i \cdot \bar D_1} , \!
  \sum\limits_{K_i \cdot \S=0} \!\! a_ie^{tK_i \cdot \bar D_1} )
  {\left(\begin{array}{ccc} {32 \over (a_3+4a_2)^2} &0&0 \\ 0&
    -{32 \over (a_3-4a_2)^2} &0 \\ 0&0&0
  \end{array} \right)}
  {\left(\begin{array}{c} \! \sum\limits_{L_j \cdot \S=2} \! \! b_je^{tL_j
      \cdot \bar D_2} \\
  \! \sum\limits_{L_j \cdot \S=-2} \! \!\! b_je^{tL_j \cdot \bar D_2}  \\
  \! \sum\limits_{L_j \cdot \S=0} \! \! b_je^{tL_j \cdot \bar D_2}
    \end{array} \right)}
$$
Here we use the fact that the third row and third column are zero.

\begin{thm}
  \label{thm:3}
  If we have $D^2 = \bar D_1^2+\bar D_2^2$, $D \cdot \S =1$ then the
  square matrix in equation~\eqref{E} is
  $$ {\left(\begin{array}{ccc} 32 e^{2t} &0&0 \\ 0&
    -32 e^{-2t} &0 \\ 0&0&0 \end{array} \right)} $$
  So $$\DD^w_X(e^{tD}) = e^{Q(tD)/2}(\sum_{K_i \cdot \S=L_j \cdot
    \S = 2} 32a_ib_j \, e^{(K_i \cdot \bar D_1 +  L_j \cdot \bar D_2
    +2)t} +\sum_{K_i \cdot
    \S=L_j \cdot \S= -2} -32a_ib_j \, e^{(K_i \cdot \bar D_1 +  L_j
    \cdot \bar D_2-2)t})$$
\end{thm}

\begin{pf}
  Here it is enough to find examples of manifolds $\bar X_1$, $\bar
  X_2$ and $X$ whose basic classes are known. Instead we use an
  indirect argument.
  Since all the manifolds involved are of simple type, the non-zero
  entries of the matrix are finite sums of exponentials, i.e.
  $$ {\left(\begin{array}{ccc} \sum c_n\, e^{nt} &0&0 \\ 0&
     \sum d_n\, e^{nt}&0 \\ 0&0&0 \end{array} \right)} $$
  Now we evaluate the series on $tD+r_1\a_1 +r_2 \a_2$, for $\a_i \in
  H^2(X_i ; \Z)$, put $t=0$ and use theorem~\ref{thm:2} to get $\sum
  c_n =32$ and $\sum d_n =-32$.
  Let $S=\CP^2 \# 10\overline{\CP}^2$ the rational elliptic surface
  blown-up once. Denote by $E_1, \dots, E_{10}$ the exceptional
  divisors and let $T_1=C-E_1- \cdots E_9$, $T_2=C-E_1- \cdots E_8-
  E_{10}$ where $C$ is the cubic curve in $\CP^2$. So $T_1$ and $T_2$
  can be represented by smooth tori of self-intersection zero and with
  $T_1 \cdot T_2 =1$. We can glue two copies of $S$ along $T_1$. The
  result is a K3 surface $S \#_{T_1} S$
  blown-up twice. The $T_2$ pieces glue
  together to give a genus $2$ Riemann surface of self-intersection
  zero $\S_2$ which intersects $T_1$ in one point.
  Now set $X= (S \#_{T_1} S ) \#_{\S_2} (S \#_{T_1} S )$, call
  $\S=\S_2$ and get $D$ piecing together both $T_1$ in $S \#_{T_1} S$.
  So (choose for instance $w=E_1$ in $S$)
  $$D_X(e^{tD+s\S}) = e^{Q(tD+s\S)/2}(\sum_{K_i \cdot \S=L_j \cdot
    \S = 2} c_na_ib_j \, e^{2s+nt} +\sum_{K_i \cdot
    \S=L_j \cdot \S= -2} d_na_ib_j \, e^{-2s+nt}) =$$
  $$ =  e^{ts} (\sum {c_n \over 16} e^{2s+nt} +\sum {d_n \over 16}
  e^{-2s+nt}) $$
  since $T_1$ evaluates $0$ on basic classes being a torus of
  self-intersection zero (the coefficient $1 \over 16$ appears from
  the explicit computation of the basic classes of the K3 surface
  blown-up in two points).
  The trick is now to use the symmetricity fact that $X= (S \#_{T_2} S
  ) \#_{\S_1}  (S \#_{T_2} S )$  where $\S_1$ comes from gluing
  together both $T_1$. Under this diffeomorphism $D=\S_1$ and
  $\S$ comes from piecing together both $T_2$ in $S \#_{T_2} S$. Hence
  $$D_X(e^{tD+s\S}) = e^{ts} (\sum {c_n \over 16} e^{2t+ns} +\sum {d_n
    \over 16}   e^{-2t+ns}) $$
  From here we deduce that $c_n=0$ unless $n=2$ and $d_n=0$ unless
  $n=-2$. Hence the result.
\end{pf}

\begin{rem}
  Note that when $\bar X_1$ is of simple type there is not summand in
  $\DD^w_{\bar X_1}$ corresponding to $s$. Therefore
  $ 0= \phi_1(a_3-16a_1) + \phi_3(a_1 - {1\over 16} a_3)$  and then
  $(\phi_3-16\phi_1)(a_3-16a_1) =0$, so $\phi_3 =16\phi_1$. In
  particular, for any $X_2$, the manifold $X=X_1 \cup_Y X_2$ is of
  simple type and, coherently, has a
  series without the coefficient corresponding to $s$.
  In this case we have
  $$ D_X(e^{tD+s\S}) = $$
  $$= {1 \over 4}
  (0 ,-16\psi_0    +\psi_2 , {1 \over 4}\psi_2 \sinh 2s + {1 \over 16}
  \psi_3 \cosh 2s, \psi_0 - { 1\over 16}\psi_2 + { 1\over 16}\psi_2
  \cosh 2s + {1 \over 64} \psi_3 \sinh 2s ) {\left( \begin{array}{c}
  \phi_0 \\ \phi_1 \\ \phi_2 \\ \phi_3 \end{array} \right)} = $$
  $$={1 \over 64}(\psi_3-16\psi_1) (16\phi_0 - \phi_2) + {1 \over 256}
  (\psi_2, \psi_3) {\left(\begin{array}{cc} 16\sinh 2s & 4 \cosh2s \\
  4\cosh2s & \sinh 2s \end{array} \right)} {\left(\begin{array}{c} \phi_2
  \\ \phi_3 \end{array} \right)}$$
  So if both $\bar X_1$ and $\bar X_2$ are of simple type, we get
  $$ \DD_X^w(tD+s\S)={1 \over 32}
  (\psi_1, \psi_2) {\left(\begin{array}{cc} 16\sinh 2s & 4 \cosh2s \\
  4\cosh2s & \sinh 2s \end{array} \right)} {\left(\begin{array}{c} \phi_1
  \\ \phi_2 \end{array} \right)}$$
\end{rem}

\begin{rem}
  \label{rem:tr}
  Suppose we do not want to use lemma~\ref{lem:PS} about the explicit
  description of $P_{\S}$. Then we could argue as follows. We have,
  instead of corollary~\ref{cor:QS},
  $$D_X(e^{s\S +tD})= <D_{X_1}^w (e^{tD_1}), Q_{\S}
  D_{X_2}^w(e^{tD_2})>$$ for some symmetric map $ Q_{\S}$.
  We go through the same steps as before without the knowledge of the
  matrices to reach the matrix $A$. For proving its invertibility it
  is enough to
  find four linearly independent vectors in the left hand side of
  formula~\eqref{C} for $\psi_i=a_i$ (i.e. $X_2=A$).
  For this we use
  \begin{itemize}
  \item $X$ a K3 surface, $\S$ a tight torus with an added trivial
    handle to make it of genus $2$. The vector we get is $(0,0,1,0)$.
  \item $X = \CP^1 \x \S$ to get a vector with non-zero last
    component. Since this manifold has $b_1 \neq 0$ we have to use the
    odd part of the Floer homology in the computations. This produces
    an extra term in the series $\DD_X^w(s\S+tD)$ of the form $f(t)$.
    The only fact that we need is that when we set $t=0$ there is a
    summand which is a multiple of $s$.
  \item $X$ a K3 surface blown-up twice, $\S = S-E_1 -E_2$ for $S$ a
    tight surface of genus $2$ in K3, $w=E_1$, $D$ a cohomology class
    coming form the K3 part such that $D \cdot S =1$. We get $(1, 1,
    {\text something} ,0)$.
  \item $X$, $\S$, $D$ as before, $w=R+ E_1$, where $R$ comes from the
    K3, has self-intersection $1$
    and is orthogonal to $D$ and $S$. We get $(1,-1, {\text something} ,0)$.
  \end{itemize}
  Having reached this point we know of the existence of a universal
  matrix as in~\eqref{D}. For the simple type case we get something like
  in~\eqref{E} with an unknown $3$x$3$-matrix. This matrix is diagonal
  since obviously it is always the case $K_i \cdot \S = L_j \cdot \S$.
  Now consider the case in which both $\bar X_i$
  and $\S_i$ are as in the first example above. Then $X = \bar X_1 \# \bar
  X_2$ splits off a ${\Bbb S}^2 \x {\Bbb S}^2$, so the invariants are
  zero. Therefore the third diagonal entry is zero and the rest of the
  argument remains intact.
\end{rem}

\subsection*{The case of genus $g >2$}

Here we propose a way of tackling conjecture~\ref{conj}.

Call $HF$ to the Floer homology of $\M$ and let $u=\mu(pt)$,
$h=\mu(\S)$ and $\G = \sum \mu(\a_{2i})\mu(\a_{2i+1})$ be
the generators of the invariant part of $HF$. Actually this invariant
part is generated as a vector space by $u^ih^j\G^p$ with $i+2p < g$
and $j+2p <g$. Now
we define $I$ to be the ideal in $HF$ generated by the image of
$H^1(\S)$ under $\mu$. The space $HF/ I$ is generated by elements
of the form $u^ih^j$ with $i <g$, $j<g$ (in principle they might not
be linearly independent).
Consider $V$ any subspace of $HF$ containing the orthogonal
complement $I^{\perp}$ of $I$ such that it has generators $e_{ij}=
u^ih^j \pmod I$, $i <g$, $j<g$. The dimension of $V$ is $N=g^2$. We
decompose $HF= V \oplus W$ with $W= V^{\perp} \subset I$.

Now we write
$$ E= e^{sh+\lambda u+ \a\G}= \sum f_{ijp}(s, \lambda, \a) u^ih^j\G^p$$
we have that for every relation $R(h,u,\G)=0$ it is $R( \pari{s},
\pari{\lambda}, \pari{\a}) E=0$ and so $R( \pari{s},
\pari{\lambda}, \pari{\a}) f_{ijp} =0$.
Note also that $ \pari{s^i}
\pari{\lambda^j} \pari{\a^p}  f_{i'j'p'} (0,0,0) =
\delta^{ijp}_{i'j'p'}$  for $i+2p < g$
and $j+2p <g$. So the $f_{ijp}$ are linearly independent functions.
$E$ defines a map from $V$ to $V$ (which we keep on calling $E$)
by multiplication followed by orthogonal projection. This map is of
the form $N_c \, g_c(s,\lambda, \a)$, where $N_c$ are constant
endomorphisms of $V$, $c =(i,j)$, $0 \le i,j < g$ and $g_c(s,
\lambda, \a)$ are linearly independent functions.

Let $X_1$ be
open manifold
with cylindrical end $Y= \Y$ and $D \in
H^2(X_1, \partial X_1 ; \Z)$ with $\partial D= \SS$. Then
$D_{X_1}(D^i) \in HF$ has component $D_{X_1}(D^i)_V$ in $V$ and
$D_{X_1}(D^i)_W$  in $W$. When $X_1$ has $b_1 =0$ one has
$<D_{X_1}(D^i), \mu(\a)>=0$ for any $1$-homology class $\a$. So
$D_{X_1}(D^i) \in V \subset HF$.
So $D_X(e^{tD})= <D_{X_1}(e^{tD_1})_V,D_{X_2}(e^{tD_2})_V>$.

Now when either of $X_i$ has $b_1=0$ one has
$D_X(e^{tD+s\S+\lambda x})= D_X(e^{tD+s\S+\lambda x+ \a \G})=
e^{ts}<D_{X_1}(e^{tD_1})_V,e_{ij}>
(<e_{ij}, e_{i'j'}>)^{-1}< e_{i'j'}, E * D_{X_2}(e^{tD_2})_V>$.
Write $\phi_a=\phi_a(t)$ for the components of $D_{X_1}(e^{tD_1})$ in
$V$ and $\psi_a=\psi_a(t)$ for the components of $D_{X_2}(e^{tD_2})$.
Then $D_X(e^{tD+s\S+\lambda x})= e^{ts}\, \phi_a(t) \, M_{abc}\,
\psi_b(t) \, g_c(s, \lambda,\a)$ for some matrices $M_{abc}$.
Now we can decompose $D_X(e^{tD+s\S+\lambda x})= e^{ts}\, D_{X,c} \,
g_c(s,\lambda, \a)$  so
$D_{X,c}= \phi_a\, (M_{abc}\, \psi_b)$ (note that the $g_c$ corresponding
to non-vanishing $D_{X,c}$ are independent of $\a$).

When $X_2=A = D^2 \x \S$, we put $a_b=\psi_b$. So we have constructed a map
\begin{eqnarray*}
  V \otimes \cal F (t)  & \to & \R^N \otimes \cal F (t) \\
  (\phi_a)_a & \mapsto & (\phi_a \, M_{abc} \, a_b)_c
\end{eqnarray*}
where $\cal F (t)$ is  the vector space of (Laurent) formal power series.
To see that $D_{\bar X_1}$ determines $\phi_a$ we need to prove the
injectivity of this linear map between vector spaces of the same
dimension.
If this were proved, the map would be an
isomorphism and hence we could mimic the argument of the case $g=2$ in
this situation to arrive to the existence of some universal $N$x$N$
matrix $P$ whose coefficients depend on $t$ and $\lambda$ such that
$$D_X(e^{tD+\lambda x})=(D_{\bar X_1}(e^{t \bar D_1})_a) (P_{ab}(t, \lambda))
(D_{\bar X_2}(e^{t \bar D_2})_b) $$

When $X$ is of simple type, $D_X(e^{tD+s\S+\lambda x})$ is a linear
combination of the functions $e^{2 \lambda}e^{(2+4n)s}$  ($-[{g \over
  2}] \leq n \leq [{g-2 \over 2}] $)
and  $e^{-2 \lambda}e^{4ns}$ ( $- [{g-1 \over 2}] \leq n \leq [{g-1
  \over 2}] $). So these functions
are among the $g_c$ (or they are
combinations of them) and without loss of generality we can suppose
they are the first $2g-1$ of the lot. With arguments as in this paper
and one non-trivial example of the gluing where the basic classes were
known, we would get the $(2g-1)$x$ (2g-1)$ minor to be (conjecturally)
$$  {\left(\begin{array}{cccc} {\text{coef }} e^{2t \pm 2 \lambda} & 0
  & \cdots  & 0 \\ 0& {\text{coef }} e^{-2t \pm 2 \lambda} & &0 \\
  \vdots & & \ddots \\   0&0& &0 \end{array} \right)} $$

We can further obtain more information from $X= \CP^1 \x \S$, but this
gives us a total of $(2g)$x$(2g)$ coefficients (very far from the
$(g^2)$x$(g^2)$ we seek for).

\section{Proof of main Theorem}
\label{sec:main}

Now we have all the information to prove theorem~\ref{thm:main}. As it was
explained in~\cite{Munoz}, we have an exact sequence
\begin{equation}
  \label{eqn:coh}
  0 \to H^1(Y;\Z) \to H^2(X; \Z) \stackrel{\pi}{\to} G \oplus H^1(\S ;\Z)
\end{equation}
where we call $G$ to the subgroup of $H^2(\bar X_1;\Z)/\Z[\S_1] \oplus
H^2(\bar X_2;\Z)/\Z[\S_2]$ consisting of elements $( \alpha_1, \alpha_2
)$ such that $\alpha_1 \cdot \S_1=\alpha_2 \cdot \S_2$ (for
interpretations of this sequence see~\cite{Munoz}).
Suppose now that both $\bar X_i$ (and hence $X$) are of simple type.
Write $\DD^w_{\bar X_1}=e^{Q/2}\sum a_i e^{K_i}$, $\DD^w_{\bar
X_2}=e^{Q/2}\sum b_j
e^{L_j}$ and $\DD^w_X=e^{Q/2}\sum c_{\kappa} e^{\kappa}$. Put $H$ for
the subgroup of $H^2(X)$ generated
by the image of $H^2(X_1) \oplus H^2(X_2)$ and a (fixed) class $D=D_1
+D_2$ such that $\bd D_1 =-\bd D_2 = \SS$ (or a non-zero multiple of
$\SS$).

First we consider classes $\bar D_i = D_i + \Delta \in H^2(\bar X_i ;
\Z)$. There is an indeterminacy as $\bar D_i $
is defined up to addition of multiples of $\S$. So we can impose $D^2
= \bar D_1{}^2 + \bar D_2{}^2$ (so $(\bar D_1,  \bar D_2)$ is defined
up to addition of multiples of $(\S, -\S)$ ). Now we glue $K_i$ and
$L_j$ (whenever $K_i \cdot \S =L_j \cdot \S= \pm (2g-2)$ ),
to get $\kappa_{ij}$.
The indeterminacy comes
this time from adding elements in $H^1(Y;\Z) $. We impose
\begin{equation}
  \kappa_{ij} \cdot D = K_i
  \cdot \bar D_1 + L_j \cdot \bar D_2 \pm 2
  \label{G}
\end{equation}
depending on whether $K_i \cdot \S =L_j \cdot \S= \pm (2g-2)$ (when
$g=1$ both possibilities must be considered). This prevents additions of
multiples of $\S$ (note that the right hand side of~\eqref{G} makes
sense). So there is still an indeterminacy coming from
addition of elements in $H^1 (\S ;\Z) \otimes H^0( \SS ;\Z) \subset
H^1(Y;\Z)$.

\begin{rem}
  The natural condition to impose on $\kappa_{ij}$ is that
  $(\kappa_{ij})^2 = K_i^2 + L_j^2 + 8(g-1)$.
  This is equivalent to the
  above whenever $g>1$.  Instead, the one
  used above is the correct one even when $g=1$.
\end{rem}

\begin{rem}
  Note that $\kappa_{ij}$ lies in $\pi^{-1}(K_i, L_j)$ by construction.
\end{rem}

For a pair $(K_i, L_j)$, all the possible $\kappa_{ij}$ restrict to
the same function $f$ in $H$. So
$$\DD^w_X |_H =e^{Q/2} \sum_{f \in H^*} (\sum_{\{ \kappa / \kappa|_H=f
  \} } c_{\kappa} ) e^f $$
Note that two cohomology classes restrict to the same $f$ if and only
if they have the same image under $\pi$ followed by projection to $G$
and they have the same pairing with $D$.

Now we return to the case of $g=2$. From what we learn in
theorems~\ref{thm:2} and~\ref{thm:3} we have
$$\DD^w_X |_H =e^{Q/2} \sum_{K_i  \cdot \S=L_j \cdot \S= 2}
  32 a_ib_j e^{\kappa_{ij}}+\sum_{K_i \cdot \S=L_j \cdot \S= -2}
  -32 a_ib_j e^{\kappa_{ij}}$$
from where we get
\begin{equation}
  \begin{cases}
   \sum\limits_{ \scriptsize{ \begin{array}{c} \kappa \in \pi^{-1}(K_i,L_j) \\
  \kappa^2 = K_i^2 +L_j ^2 +8 \end{array}}}
   c_{\kappa} = \pm 32a_ib_j& \text {if
  $K_i \cdot \S =L_j \cdot \S = \pm 2$} \\
    \sum\limits_{ \kappa \in \pi^{-1}(K_i,L_j)} c_{\kappa} = 0 & \text {if
  $K_i \cdot \S =L_j \cdot \S = 0$ } \\
    \sum\limits_{ \kappa \in \pi^{-1}(K,L)} c_{\kappa} = 0 & \text {if
  $(K,L) \neq (K_i, L_j)$}
   \end{cases}
  \label{F}
\end{equation}
Note that for any basic class $\kappa$ one has $\kappa \cdot
T_{\gamma} =0$ (where $T_{\gamma} = \SS \x \gamma \subset \Y$
  for $\gamma \in H_1(\S; \Z)$). Therefore $\pi(\kappa)$ lies in $G
  \subset G \oplus H^1(\S; \Z)$.
Also note that the condition $\kappa^2 = K_i^2 +L_j ^2 +8$ in the
first case means that the possible $\kappa$ in the sum differ by
addition of elements in $H^1 (\S ;\Z) \otimes H^0( \SS ;\Z) \subset
H^1(Y;\Z)$.

\begin{cor}
  Let $\bar X_i$ as before and $g=2$.
  Suppose that for every cycle $\gamma \in H^1(\S ; \Z)$ there exists
  a $(-1)$-embedded disc in both $X_i$ bounding $\gamma$, then the
  basic classes $\kappa$ of $X$
  are in one-to-one correspondence with
  pairs of basic classes $(\kappa_1, \kappa_2)$ for $\bar X_1$ and $\bar X_2$
  respectively, such that $\kappa_1 \cdot \S_1 =\kappa_2 \cdot \S_2
  =\pm 2$. Moreover, $\kappa$ is determined in an explicit way.
\end{cor}

\begin{pf}
  As explained in~\cite{Munoz}, there is a splitting $H^2(X ;\Z)=V
  \oplus V^{\perp}$ where $V$ is generated by tori
  $T_{\gamma} = \SS \x \gamma \subset \Y$
  (for $\gamma \in H_1(\S; \Z)$) and spheres $D_{\beta}$ of
  self-intersection $-2$ (with the property  $D_{\beta} \cap Y = \beta$). There
  is an exact sequence
  $$0 \ar \Z[\S] \ar V^{\perp} \stackrel{\pi}{\ar} G$$
  Now let $\kappa$ be a basic class for $X$. We argue as in~\cite{KM} that
  $\kappa \in V^{\perp}$. Thus in the summation of~\eqref{F} only one term
  is non-vanishing when $K_i \cdot \S =L_j \cdot \S = \pm 2$.
  $\kappa$ is characterised as the only class orthogonal to all
  $T_{\gamma}$  and $D_{\beta}$ such that $\kappa^2 = K_i^2 +L_j ^2
  \pm 8$.
  When $K_i \cdot \S =L_j \cdot \S =0$ there is only one term in the
  summation since the condition~\eqref{G}. Therefore there are no basic
  classes such that $\kappa \cdot \S = 0$.
\end{pf}

\end{document}